\documentclass[12pt]{article}
\usepackage{amsmath}
\usepackage{amscd}
\usepackage{amsfonts}
\usepackage{epsfig}
\usepackage{theorem}

\newtheorem{theorem}{Theorem}[section]

\newtheorem{lemma}{Lemma}[section]

\newtheorem{definition}{Definition}[section]

\newtheorem{conjecture}{Conjecture}[section]

\def \proof {\vspace{-.4cm}{\bf Proof \quad}}

\def \epsilon {\varepsilon}

\def \phi {\varphi}

\def \G {{\cal G}}

\def \Gm {{\cal G}_m}

\def \J {{\cal J}}

\def \Jm {{\cal J}_m}

\def \Jn {{\cal J}_n}

\def \Jmt {\tilde{\cal J}_m}

\def \F {{\cal F}}

\def \Fm {{\cal F}_m}

\def \Fn {{\cal F}_n}

\def \K {{\cal K}}

\def \Km {{\cal K}_m}

\def \Kmt {\tilde{\cal K}_m}

\def \A {{\cal A}_\infty}

\def \Am {{\cal A}_{\infty, m}}

\def \S {{\mathcal S}}

\def \Smn {{\mathcal S}_m^n}

\def \Sn {{\mathcal S}^n}

\def \G {{\cal G}}

\def \calG {{\mathcal G}}

\def \C {\mathbb C}

\def \cbar {{\overline {\mathbb C}}}

\def \N {\mathbb N}

\def \Q {\mathbb Q}

\def \nocbar {{\mathbb N}_0 \times {\overline {\mathbb C}}}

\def \G {{\cal G}}

\def \Pm {\{P_m \}_{m=1}^\infty}

\def \Pmt {\{\tilde P_m \}_{m=1}^\infty}

\def \linfty {l_\infty}

\def \meas {{\mathrm m}}

\def \dzbardz {\frac{{\mathrm d}\overline z}{{\mathrm d}z}}

\title{Non-Autonomous Julia Sets with Invariant Sequences of Measurable Line Fields}

\author{Mark Comerford\\ 
Department of Mathematics\\
University of Rhode Island\\
5 Lippitt Road, Room 102F\\ 
RI 02881, USA\\ \\
email: {\tt mcomerford@math.uri.edu}}

\begin{document}  

\maketitle

\begin{abstract}The no invariant line fields conjecture is one of the main outstanding problems in traditional complex dynamics. In this paper we consider non-autonomous iteration where one works with compositions of sequences of polynomials with suitable bounds on the degrees and coefficients. We show that the natural generalization of the no invariant line fields conjecture to this setting is not true. In particular, we construct a sequence of quadratic polynomials whose iterated Julia sets all have positive area and which has an invariant  sequence of measurable line fields whose supports are these iterated Julia sets with at most countably many points removed. 
\end{abstract}

\section{Introduction}

The no invariant line fields conjecture is one of the most important unsolved problems in complex dynamics, especially because of its relation to the question of density of hyperbolicity (e.g. \cite{MSS}). 
In this paper we exhibit a counterexample which shows that the obvious generalization of the no invariant line fields conjecture to non-autonomous polynomial iteration is false. In non-autonomous iteration, one considers composition sequences of polynomials with suitable bounds on the degrees and coefficients. About ten years ago, Peter Jones together with Michael Benedicks, Mattias Jonsson and Misha Yampolsky found an example of a sequence of quadratic polynomials whose Julia set has positive area (sketches of their proof may be found in \cite{Com1, Com2}). In this paper, we modify and considerably extend their ideas to show how one can construct a sequence of quadratic polynomials for which there is an invariant sequence of measurable line fields which are supported on very nearly the whole of the iterated Julia sets. 

We begin with the basic definitions we need in order to state the main result of this paper. In the next section we prove this theorem, together with some supporting results and make a few concluding remarks. 

\subsection{Bounded Polynomial Sequences}

Let $d \ge 2$, $M \ge 0$, $K \ge 1$ and let $\Pm$ be a sequence of polynomials where each $P_m(z) = a_{d_m,m}z^{d_m} + a_{d_m-1,m}z^{d_m-1} + \cdots \cdots + 
a_{1,m}z + a_{0,m}$ is a polynomial of degree $2 \le d_m \le d$ whose coefficients satisfy
\[\qquad 1/K \le |a_{d_m,m}| < K,\quad m \ge 1, \:\:\quad \quad |a_{k,m}| \le M,\quad m \ge 1,\: 0 \le k \le d_m -1. \]
Such sequences are called \emph{bounded sequences of polynomials} or simply \emph{bounded sequences}. Such a sequence can be thought of as a skew product on $\cbar$ over the non-negative integers $\N_0$, or equivalently as a mapping $F$ of the set $\nocbar$ to itself, given by $F(m, z) := (m+1, P_{m+1}(z))$.

For each $1 \le m$, let $Q_m$ be the composition $P_m \circ \cdots \cdots \circ P_2 \circ P_1$ and for each $0 \le m \le n$, let $Q_{m,n}$ be the composition $P_n \circ \cdots \cdots \circ P_{m+2} \circ P_{m+1}$ (where we let each $Q_{m,m}$ be the identity). Let the degrees of these compositions be $D_m$ and $D_{m,n}$ respectively so that $D_m = \prod_{i=1}^m d_i$, $D_{m,n} = \prod_{i=m+1}^n d_i$.

For each $m \ge 0$ define the \emph{$m$th iterated Fatou set} $\Fm$ by 
\[ \Fm = \{z \in \cbar : \{Q_{m,n}\}_{n=m}^\infty \;
\mbox{is a normal family on some neighbourhood of}\; z \}\]
where we take our neighbourhoods with respect to the spherical topology on $\cbar$ 
and let the \emph{$m$th iterated Julia set} $\Jm$ be the complement $\cbar \setminus \Fm$. At time $m =0$, we call the corresponding iterated Fatou and Julia sets simply the \emph{Fatou} and \emph{Julia sets} for our sequence and designate them by $\F$ and $\J$ respectively.

One can easily show that the iterated Fatou and Julia sets are completely invariant in the following sense. 

\begin{theorem}
For each $0 \le m \le n$, $Q_{m,n}(\Fm) = \Fn$ and $Q_{m,n}(\Jm) = \Jn$ with components of $\Fm$ being mapped surjectively onto components of $\Fn$.
\end{theorem}

This definition also has the advantage that we can find some radius $R$ 
depending only on the bounds $d$, $K$, $M$ above so that for any sequence $\Pm$ as above and any $m \ge 0$, it is easy to see that 
\[ |Q_{m,n}(z)| \to \infty \qquad \mbox{as} \quad n \to \infty, \quad |z| > R\]
which shows in particular that as for classical polynomial Julia sets, for each $m \ge 0$ there 
will be an iterated basin at infinity $\Am$ on which all points escape to infinity under
iteration. Such a radius will be called an {\it escape radius} for the 
coefficient bound $M$. Note that the maximum principle shows that, just as in 
the classical case, there can be only one component on which $\infty$ is 
a normal limit function and so the sets $\Am$ are 
completely invariant in the sense given in Theorem 1.1 above.

The complement of $\Am$ is called the {\it iterated filled Julia set at time $m$} for
the sequence $\Pm$ and is denoted by $\Km$. As above, when $m = 0$, we refer simply to the basin at infinity and filled Julia set for our sequence and denote them by $\A$ and $\K$ respectively. 

The iterated filled Julia sets are then also clearly completely invariant in the sense given above. Also, the same 
argument using Montel's theorem as in the classical case shows that, for each $m \ge 0$, 
$\partial \Km = \Jm$ as one might expect.  If $R$ is an escape radius as above, for $0 \le m \le n$, let us call the set $Q_{m,n}^{\circ -1}({\mathrm D}(0,R))$ the \emph{$n$th survival set at time $m$} and denote it by $\Smn$. Clearly, the sets $\Smn$ are decreasing in $n$ and $\Km = \bigcap_{n \ge m} \Smn$ for each $m \ge 0$. Finally, following our earlier notation, at time $0$, we denote the sets $\S_0^n$ simply by $\Sn$. 

Another important example of complete invariance is the grand orbit of a set. 
For a bounded sequence of polynomials 
$\Pm$, and some set $X \subset \cbar$ which appears at time $m \ge 0$, for $i \ge 0$ we obtain the {\it grand orbit at time $i$} of $X$, $\calG_i (X)$ given by 
\[\calG_i (X) := \{Y \subset \cbar : Q_{i,j}(Y) = Q_{m,j}(X) \quad \mbox{for some} \quad j \ge \max (m,i)\}.\]
If $i=m$, we call $\calG_m(X)$ the {\it immediate grand orbit} of $U$ (at time $m$). Lastly, if $m = 0$, we denote the grand orbit of $X$ at time $0$ simply by $\G(X)$.

\subsection{Invariant Line Fields}

Recall that a line field supported on a subset $E$ of $\cbar$ is choice of a real line through the origin in the tangent space ${\mathrm T}_e(z)$ for each $e \in E$. Equivalently (e.g. in the book of McMullen \cite{McM}), it is a Beltrami form $\mu = \mu(z) \dzbardz$ where $|\mu(z)| = 1$ everywhere on $E$. 
Following McMullen, we say $\mu$ is invariant under a rational map $R$ if $R^* \mu = \mu$ almost everywhere, $| \mu(z) | = 1$ on a set of positive measure and $\mu(z)$ vanishes elsewhere. We now state the no invariant line fields conjecture. For details on the statement, the reader is referred to \cite{McM}.

\begin{conjecture}
A rational map $R$ carries no invariant line field on its Julia set, except when $R$ is double covered by an integral torus endomorphism. 
\end{conjecture}

Some progress has been made on this; for example by Wang for the case of meromorphic invariant line fields \cite{Wang} and by Rempe and Van Strien for transcendental meromorphic functions with certain restrictions \cite{RVS}. On the other hand, Eremenko and Lyubich show in \cite{EL} that there are certain entire transcendental functions which have invariant line fields supported on their Julia sets. However, the conjecture remains open, even for quadratic polynomials of the form $z^2 + c$ (where it is equivalent to density of hyperbolicity). 

We now proceed to define the analogue of an invariant line field for non-autonomous polynomial iteration. 

\begin{definition}
Let $\{\mu_m\}_{m=0}^\infty$ be a sequence of line fields and let $\Pm$ be a bounded sequence of polynomials. We say such a sequence is \emph{invariant} under $\Pm$ if for every $0 \le m < n$, we have $Q_{m,n}^*(\mu_n)(z) = \mu_m(z)$. 
\end{definition}

The main result of this paper is then as follows. 

\begin{theorem}
There exists a bounded sequence $\Pm$ of quadratic polynomials for which we have an invariant sequence $\{\mu_m\}_{m=0}^\infty$ of measurable line fields where for each $m \ge 0$, the support of $\mu_m$ is a subset $\Jmt$ of $\Jm$ where $\Jm \setminus \Jmt$ is at most countable. 
\end{theorem}

The author wishes to gratefully acknowledge the inspiration and advice of Lasse Rempe who first posed the question on invariant line fields and Hiroki Sumi who gave the correct definition for the non-autonomous version, which together led to the writing of this paper.

\section{Proof of the Main Result}

The main obstacle to constructing a non-zero measurable invariant line field is lack of injectivity as this makes it difficult to push forward the Beltrami differential. The basic idea is to start with a disc and construct the sequence so that there is a dense subset on which the iterates escape to infinity and whose measure is strictly less than that of the original disc. All points on the residual set will have bounded orbits and as this set contains no interior, it must be a subset of the Julia set. If the iterates on this part of the Julia set are injective, then we will eventually be able to construct an invariant sequence of measurable line fields whose supports are contained in the iterated Julia sets. 

The proof is based on the dynamics associated with two polynomials: $P^1(z) = \tfrac{1}{2}z (1 + z)$ and 
$P^2(z)  = z(1 + z)$. $P^1$ has an attracting fixed point at $0$ and a simple calculation shows that $|P^1(z)| < |z|$ provided $|z| < 1$. Since the critical point of $P^1$ is at $-\tfrac{1}{2}$, we may extend the Schr\"oder linearization conjugacy to all of ${\mathrm D}(0, \tfrac{1}{2})$ (e.g. \cite{CG}). It follows that the iterates of $P^1$ will then be injective on the smaller disc ${\mathrm D}(0, \tfrac{1}{3})$.

$P^1$ is hyperbolic and has just one bounded Fatou component which is the basin of attraction of the fixed point at $0$.  If we let $U^1$ denote this component and $\K^1 = \overline {U^1}$ the corresponding filled Julia set, then ${\mathrm D}(0, \tfrac{1}{3}) \subset U^1$ and since ${(P^1)}^{\circ -1}({\mathrm D}(0, \tfrac{1}{3}))$ compactly contains ${\mathrm D}(0, \tfrac{1}{3})$, it follows that the sets ${(P^1)}^{\circ -n}({\mathrm D}(0, \tfrac{1}{3}))$ are increasing in $n$ and expand to fill out all of $U^1$ as $n$ tends to infinity. Using this and the hyperbolicity of $P^1$, we obtain the following. Although this result is simple and its proof straightforward, it is the key to ensuring that a sequence of non-vanishing line fields can be defined on almost the whole of the iterated Julia sets. For a measurable subset $X$ of $\C$, we use the standard notation $\meas(X)$ to denote the two-dimensional Lebesgue measure of $X$

\begin{lemma}Let $\Pm$ be a bounded polynomial sequence with bounds $d, K,M$ as above where $K \ge 2$ and $M \ge \tfrac{1}{2}$. Then there exist $\lambda > 1$ depending only on $P^1$ and $c>0$ depending on $d,K,M$ and $P^1$ such that if $N \ge1$ is such that the first $N$ members of our sequence are all $P^1$ and  $\S^N$ as usual denotes the $N$th survival set at time $0$, then 
\vspace{.3cm}
\[ m(\S^N \setminus \K^1) < c \lambda^{-N}.\]

\end{lemma}

An important fact to note in connection with this result is that we do not need to know in advance the particular sequence of polynomials involved, only the bounds on the coefficients. This means that we can get estimates on the measure of parts of the filled Julia set, even before the sequence has been fully constructed!

Passing now to $P^2$, this polynomial has a parabolic fixed point with multiplier $1$ at $0$ and the associated Julia set is the 
well-known `cauliflower'. The repelling direction is the positive real axis while the attracting direction is the negative real axis. There is also an attracting petal $U^2$ associated with $P^2$ which for the sake of convenience we shall assume is open. $U^2$ is forward invariant and the iterates on $U^2$ are conjugate the the mapping $z \mapsto z + 1$.
Since this mapping is injective, $P^2$ is also injective on $U^2$
and so in particular the critical point of $P^2$ at $-\tfrac{1}{2}$ as well as all its preimages under $P^2$ lie outside $U^2$. 

$\C \setminus U^2$ has a cusp at $0$ in the direction of the positive real axis so that the angle subtended by the intersection of a small disc ${\mathrm D}(0,r)$ about $0$ with this cusp can be made arbitrarily small by taking the radius $r$ of this disc small enough. 
Using this we can prove the following simple lemma. Roughly speaking, this result says that if we apply a high iterate of $P^1$ to an open subset of ${\mathrm D}(0, \tfrac{1}{3})$ and then carry out a translation so that the image set still meets the tip of the cusp of $\C \setminus U^2$ at $0$, then the measure of the preimage of the intersection with this cusp can be made arbitrarily small. For $X \subset \C$ and $z_0 \in \C$ fixed, let us denote the translated set $\{z \in \C: z-z_0 \in X\}$ by $X + z_0$.

\begin{lemma}Let $V$ be an open subset of ${\mathrm D}(0, \tfrac{1}{3})$ and let $\epsilon >0$. Then we can find $m_0$ large enough so that for any $m \ge m_0$ and any $z_0 \in {(P^1)}^{\circ m}(V)$
\[ \meas(V \cap {(P^1)}^{\circ -m}((\C \setminus U^2) + z_0))) < \epsilon.\]

\end{lemma}

\vspace{.2cm}
\proof ${\mathrm D}(0, \tfrac{1}{3})$ is contained within finite hyperbolic distance of $0$ in ${\mathrm D}(0, \tfrac{1}{2})$ and the iterates of $P^1$ on this larger disc are univalent. We may then apply the distortion theorems for univalent mappings (e.g. Theorem 1.6 on Page 5 of \cite{CG}) to find $k >0$ such that for any $m \ge 1$, if $r >0$ is as small as possible so that ${(P^1)}^{\circ m}(V) \subset {\mathrm D}(0,r)$, then $\meas ({(P^1)}^{\circ m}(V)) \ge kr$. Thus, given $\delta >0$, we can find $m_0$ large enough to ensure that, for any $m \ge m_0$ and any $z_0 \in V$, 
\vspace{.2cm}
\[\meas ({(P^1)}^{\circ m}(V) \cap ((\C \setminus U^2) + z_0)) < \delta \meas ({(P^1)}^{\circ m}(V)).\] 

Again as ${\mathrm D}(0, \tfrac{1}{3})$ is contained within finite hyperbolic distance of $0$ in ${\mathrm D}(0, \tfrac{1}{2})$ and the iterates of $P^1$ on ${\mathrm D}(0, \tfrac{1}{2})$ are univalent, the conclusion follows by applying the distortion theorems to the suitable inverse branches of $P^1$. $\Box$

All polynomials in our sequence $\Pm$ will be either $P^1$, $P^2$, $P^1+ c_m$ where $|c_m| < \tfrac{1}{3}$ or $\tfrac{1}{12}P_2$, the last having the property that it maps the entire filled Julia set for $P^2$ inside the disc ${\mathrm D}(0, \tfrac{1}{3})$. This clearly gives us a bounded sequence and some simple calculations show that if $|z| > 13$, then the iterates of $z$ under any sequence whose members are chosen from the above must escape to infinity. We now proceed with the proof of the main result.

{\bf Proof of Theorem 1.2 \hspace{.4cm}} As is standard in these situations, the construction of the required sequence will be an inductive one consisting of infinitely many stages of iteration, where the end of each stage sets things up for the next. 

Let  $D$ be the disc ${\mathrm D}(0, \tfrac{1}{3})$. We introduce this disc at various times in our iterative procedure, look at the inverse images at time $0$ and make use of Lemma 2.1 to say that they cover almost the entire filled Julia set.  

We also introduce sequences $\{s_n\}_{n=1}^\infty$, $\{t_n\}_{n=1}^\infty$, $\{u_n\}_{n=1}^\infty$ of natural numbers, $\{\epsilon_n\}_{n=1}^\infty$ of positive real numbers all of which are less than $\meas(D) = \tfrac{\pi}{9}$ and $\{c^n\}_{n=1}^\infty$ of complex numbers. All five of these sequences are determined as we proceed with the construction of our sequence of polynomials.  

{\bf Induction: Stage 1\hspace{.4cm}} We iterate with $P^1$  $s_1$ times, introduce the disc $D^1 = D$ at time $s_1$, iterate with $P^1$ a further $t_1$ times, apply a translation by $c^1$, then iterate with $P^2$ $u_1$ times and finally apply a dilation by $\tfrac{1}{12}$. Letting $M_1 = m_1 = s_1 + t_1 + u_1$, we have defined the first $M_1$ polynomials according to 

\[P_m = \left \{ \begin{array}{r@{\: ,\quad}l} 
P^1 & 1 \le m \le s_1\\
P^1 & s_1 < m < s_1 + t_1\\
P^1 + c^1 & m = s_1 + t_1\\
P^2 & s_1 + t_1 < m < M_1\\
\frac{1}{12}P_2 & m = M_1\\
\end{array} \right .\]

Using Lemma 2.1, we can choose $s_1$ large enough so that if we let $E^1 = {(Q_{s_1})}^{\circ -1}(D^1) = {(P^1)}^{\circ -s_1}(D^1)$, then  $\meas(\S^{s_1} \setminus E^1) < \tfrac{1}{4}$ (recall that the lemma requires us only to know in advance the bounds for our sequence of polynomials, not the actual sequence itself). Since $s_1$ is now fixed, we can choose $\epsilon_1$ small enough so that, for any subset $X$ of ${\mathrm D}(0,13)$ at time $s_1$ with $\meas(X) \le \epsilon_1$, $\meas({(Q_{s_1})}^{\circ -1}(X)) < \tfrac{1}{4}$. It then follows that if $Y$ is a subset of $D^1$ of measure at least $\tfrac{\pi}{9} - \epsilon_1$, then we must have $\meas(\S^{s_1} \setminus {(Q_{s_1})}^{\circ -1}(Y)) < \tfrac{1}{2}$.

Let $\{z^n\}_{n=1}^\infty$ be a denumeration of the countable dense subset $(\Q \times \Q) \cap D^1$ at time $s_1$ (this is the only time at which we need to introduce such a dense subset).
Using Lemma 2.2, we can choose $t_1$ large enough so that for any $z_0 \in {(P^1)}^{\circ t_1}(D^1)$, we have $\meas(D^1 \cap {(P^1)}^{\circ -t_1}((\C \setminus U^2) + z_0))) < \tfrac{\epsilon_1}{4}$. Now let $c^1 = - {(P^1)}^{\circ t_1}(z^1)$ so that the point $z^1$ is iterated at time $s_1 + t_1$ by $Q_{s_1, s_1 + t_1}$ to $0$ which is the tip of the cusp of $\C \setminus U^2$. 

For convenience of notation later on, we denote $D^1$ by $B^{1,1}$ and then 
denote by $B^{1,2}$ the subset of those points of $D^1$ whose iterates are mapped at time $s_1 + t_1$ inside the petal $U^2$ for $P^2$ and thus inside ${\mathrm D}(0, \tfrac{1}{3})$ at time $M_1$. From above we then have that $\meas(B^{1,2}) > \tfrac{\pi}{9} - \tfrac{\epsilon_1}{4}$.

Since $s_1$, $t_1$ and $c^1$ are now fixed, we can choose $u_1$ large enough so that there is a point $w^1$ at time $s_1$ within distance $1$ of $z^1$ for which $|Q_{M_1}(w^1)| > 13$ so that this point is guaranteed to escape to infinity, regardless of how we choose the subsequent polynomials for our sequence from among the four possibilities. 

For compatibility with what comes later, we note that since $D^{1,1}$ has diameter less than $1$, all points of $\{z^n\}$ are within distance $1$ of $w^1$. In particular, if $z^i$ is any one of the points of our countable dense subset above and $z^i \in B^{1,1} \setminus B^{1,2}$, then $z^i$ is clearly within distance $1$ of a point $w$ which is mapped outside ${\mathrm D}(0, 13)$ at time $M_1$ and which therefore lies in ${\mathcal A}_{\infty, s_1}$, the basin at infinity at time $s_1$ for our sequence. We then define the subset $I_1$ of $\N$ by
\[ I_1 = \{ i : z^i \in B^{1,1} \setminus B^{1,2}\}.\]
The dilation by $\tfrac{1}{12}$ which we performed at the end maps the petal $U^2$ into ${\mathrm D}(0, \tfrac{1}{3})$ which allows us to continue to the next stage of the construction. Lastly, in order to be able to construct our invariant line field at the end of the proof, we observe that the composition $Q_{s_1, M_1 }$ is injective when restricted to the set $B^{1,2}$.

{\bf Induction Hypothesis: Stage $n$\hspace{.4cm}}

We assume now that the first $n$ stages of the construction have been carried out. For each $1 \le i \le n$, stage $i$ begins as above with an iterate of $P^1$ followed by a 
further iterate with $P^1$ and then a translation. This is then followed by an iterate of $P^2$ and finally a dilation by $\tfrac{1}{12}$. For each $1 \le i \le n$, let $s_i$, $t_i$ be the number of iterations of $P^1$ and $u_i$ the number of iterations of $P^2$. Next set 
$m_i = s_i + t_i + u_i$ and $M_i = \sum_{j=1}^i m_j$ (where for convenience we set $M_0 = 0$) so that we have now constructed the first $M_n$ polynomials of our sequence. 

At each time $M_{i-1} + s_i$, we introduce copies $D^i$ of $D = {\mathrm D}(0, \tfrac{1}{3})$ and let $E^i$ be the inverse image of $D^i$ at time $0$ under the composition $Q_{M_{i-1} + s_i}$. For each $1 < j \le n - i + 2$, let $B^{i,j}$ be the set of points in $D^i$ which are mapped inside $U^2$ at time $M_{k-1} + s_k + t_k$ for each $i \le k \le i + j - 2$ (where for convenience we set $B^{i,1} = D^i$).  We also introduce subsets $I^i$, $1 \le i < n$ of $\N$ defined by 
\[ I_i = \{ k : z^k \in B^{1,i} \setminus B^{1,i+1} \}.\]
Our induction hypotheses are then as follows:

\begin{enumerate}
\item $\meas(\S^{M_{n-1} + s_n} \setminus E^n) < 2^{-n - 1}.$

\item $\epsilon_n$ is chosen such that for any measurable subset $X$ of $D = {\mathrm D}(0, 13)$ at time $M_{n-1} + s_n$ of measure $< \epsilon_n$, $\meas(Q_{M_{n-1} + s_n}^{\circ -1}(X)) < 2^{-n - 1}.$

\item $\meas(B^{i,j}) > \tfrac{\pi}{9} - \epsilon_i \sum_{k=2}^j \tfrac{1}{2^k}$, $\:2 \le j \le n -i + 2$.

\item The sets $I_i$, $1 \le i \le n$ are pairwise disjoint and\\ $\{1, 2, \ldots \ldots, n\} \subset \bigcup_{i=1}^n I_i$.

\item For each $k \subset I_i$, there is a point in ${\mathcal A}_{\infty, s_1}$ within distance $< \tfrac{1}{i}$ of $z^k$. 

\item $B^{i,j+1} \subset B^{i,j}$ for each $1 \le j \le n - i + 1$. 

\item $Q_{M_{i-1}+s_i, M_i + s_{i+1}}(B^{i,j}) \subset B^{i+1,j-1}$ for each $2 \le j \le n - i + 2$.

\item $Q_{M_{i-1} + s_i, M_j}$ is injective on the set $B^{i,j - i +2}$ for each $i \le j \le n$. 

\end{enumerate}

{\bf Induction Step: Stage $n + 1$\hspace{.4cm}}

Let $s_{n+1}, t_{n+1}$, $u_{n+1}$ be natural numbers, $\epsilon_{n+1} >0$ and $c_{n+1} \in \C$, all of which we will choose later. We iterate with $P^1$  $s_{n+1}$ times, iterate again with $P^1$ a further $t_{n+1}$ times and apply a translation by $c^{n+1}$. We then iterate with $P^2$ $u_{n+1}$ times and apply a dilation by $\tfrac{1}{12}$. Now let $m_{n+1} = s_{n+1} + t_{n+1} + u_{n+1}$, and $M_{n+1} = M_n + m_{n+1} = \sum_{i=1}^{n+1} m_i$. We have then defined $P_m$ for $M_n +1 \le m \le M_{n+1}$ by

\[P_m = \left \{ \begin{array}{r@{\: ,\quad}l} 
P^1 & M_n + 1 \le m \le M_n + s_{n+1}\\
P^1 & M_n + s_{n+1} < m < M_n + s_{n+1} + t_{n+1}\\
P^1 + c^{n+1} & m = M_n + s_{n+1} + t_{n+1}\\
P^2 & M_n + s_{n+1} + t_{n+1} < m < M_{n+1}\\
\frac{1}{12}P^2 & m = M_{n+1}\\
\end{array} \right .\]

Now let $D^{n+1} = D$ be introduced at time $M_n + s_{n+1}$. Since the first $M_n$ polynomials of our sequence are already fixed, using Lemma 2.1, we can make $s_{n+1}$ sufficiently large so that if we let $E^{n+1}$ be the inverse image of $D^{n+1}$ under $Q_{M_n + s_{n+1}}$, then $\meas(\S^{M_n + s_{n+1}} \setminus E^{n+1}) < 2^{-n - 2}$ which satisfies the first of our induction hypotheses for $n+1$.

We have now chosen the first $M_n + s_{n+1}$ polynomials of our sequence and so we may pick $\epsilon_{n+1}$ sufficiently small so that for any measurable subset $X$ of $D = {\mathrm D}(0, 13)$ at time $M_n + s_{n+1}$ of measure $\le \epsilon_{n+1}$, 
$\meas(Q_{M_{n} + s_{n+1}}^{\circ -1}(X)) < 2^{-n - 2}$ and so we have satisfied our second induction hypothesis for $n+1$.

Set $B^{1,n+1} = D^{n+1}$ and for each $1 \le i \le n+1$, set $B^{i, n - i + 3}$ to be the set of points in $D^i$ which are mapped inside $U^2$ at time $M_{k-1} + s_k + t_k$ for each $i \le k \le n +1$ (note that these sets are only fully defined once we have chosen 
the integer $t_{n+1}$ and the complex number $c^{n+1}$).  

Again as the first $M_n + s_{n+1}$ polynomials are already fixed, we may use Lemma 2.2 to choose $t_{n+1}$ large enough so that for any $0 \le i \le n+1$, if we let $\tilde Q_{M_{i-1} + s_i, M_n + s_{n+1} + t_{n+1}}$ denote the polynomial composition ${(P^1)}^{\circ t_{n+1}} \circ Q_{M_{i-1} + s_i, M_n + s_{n+1}}$, then for any point $z_0 \in \tilde Q_{M_{i-1} + s_i, M_n + s_{n+1} + t_{n+1}}(B^{i,n-i+3}))$, we have 

\[ \meas\left (B^{i,n-i +3} \cap {\tilde Q_{M_{i-1} + s_i, M_n + s_{n+1} + t_{n+1}}}^{\circ -1}((\C \setminus U^2) + z_0))\right ) < \frac{\epsilon_i}{2^{n-i+3}}.\] 

By the third induction hypothesis for $n$, $B^{1,n+2}$ is an open set of positive area whence $\bigcup_{1 \le i \le n} I_i \ne \N$. We can thus 
let $j_{n+1}$ be the smallest natural number such that $z^{j_{n+1}} \notin \cup_{i = 1}^n I_i$. Note that by the fourth induction hypothesis $j_{n+1} \ge n+1$. Now let the constant $c^{n+1}$ be
$ - \tilde Q_{s_1, M_n + s_{n+1} + t_{n+1}}(z^{j_{n+1}}))$ so that the point $z^{j_{n+1}}$ is iterated at time $M_n + s_{n+1} + t_{n+1}$ to the tip of the cusp of $\C \setminus U^2$ at $0$. 

The above estimate on the measure of the inverse images of translates of the set $\C \setminus U^2$ and the third induction hypothesis for $n$ ensure this induction hypothesis is also satisfied for $n+1$. Now let $I_{n+1}$ be the set of natural numbers $\{k: z^k \in B^{1,n+1} \setminus B^{1,n + 2} \}$. The sets $I_1, I_2, \ldots \ldots, I_n, I_{n+1}$ are then disjoint which satisfies the first part of the fourth induction hypothesis for $n+1$. Also, the second part of the fourth induction hypothesis for $n$ and the fact that $j_{n+1} \ge n+1$ give us the second part of this hypothesis for $n+1$.

Any point of $B^{1, n+1} \setminus B^{1, n+2}$ is mapped into the intersection of a small disc ${\mathrm D}(0,r)$ with the complement of the petal $U^2$. Although $U^2$ does not meet the Julia set for $P^2$ except at $0$, as $\C \setminus U^2$ meets $0$ in a cusp, by taking $t_{n+1}$ and then $u_{n+1}$ large enough, it follows that any point in ${\mathrm D}(0,r) \setminus U^2$ can be made arbitrarily close to points which escape outside ${\mathrm D}(0,13)$. Again as $\C \setminus U^2$ meets $0$ in a cusp, the distance of points in ${\mathrm D}(0,r) \setminus U^2$ to points which escape to infinity can be made arbitrarily small, not just in absolute terms, but also relative to $r$ (note that we may have to make $t_{n+1}$ larger here, but this will not affect the fact that the third and fourth hypotheses hold for $n+1$).

As in the proof of Lemma 2.2, ${\mathrm D}(0, \tfrac{1}{3})$ is contained within finite hyperbolic distance of $0$ in ${\mathrm D}(0, \tfrac{1}{2})$ and $Q_{M_n+s_{n+1}, M_n + s_{n+1} + t_{n+1}} = {(P^1)}^{\circ t_{n+1}} + c^{n+1}$ is univalent on the larger disc. We may then apply the distortion theorems for univalent mappings (e.g. Theorem 1.6 on Page 3 of \cite{CG}), to the appropriate inverse branch of $Q_{M_n + s_{n+1}, M_n + s_{n+1} + t_{n+1}}$ on $Q_{M_n+s_{n+1}, M_n + s_{n+1} + t_{n+1}}({\mathrm D}(0, \tfrac{1}{2}))$. Since $Q_{M_n + s_{n+1}}$ is already fixed, it follows that for any member $z^j$ of our sequence $\{z^n\}$ which lies in $B^{1, n+1} \setminus B^{1, n+2}$, we can find a point of ${\mathcal A}_{\infty, s_1}$ within distance less than $\tfrac{1}{n+1}$ of $z^j$. Thus, we have satisfied the fifth induction hypothesis for $n+1$.

The sixth induction hypothesis for $n+1$ follows easily from the case for $n$ and the definition of the new sets $B^{i,n-i+3}$. 
The last polynomial $P_{M_n}$ from stage $n$ maps $U^2$ into a compact subset of ${\mathrm D}(0,\tfrac{1}{3})$ and the seventh induction hypothesis follow easily from this and the corresponding case for $n$. The eighth induction hypothesis for $n+1$ also follows from this and the corresponding case for $n$ together with the fact that $Q_{M_n, M_{n+1}}$ is clearly injective on ${\mathrm D}(0,\tfrac{1}{3})$. 

This completes the induction part of the argument. We now show that the Julia set $\J_{s_1}$ at time $s_1$ has positive Lebesgue measure after which we will proceed to construct the desired sequence of line fields. 

By part 6. of the induction statement, for each $n \ge 1$, the sets $B^{n,i}$, $i \ge 1$ are nested and measurable (in fact, they are open). On setting $B^n = \bigcap_{i \ge 1} B^{n,i}$, by part 3. of the induction statement, we have that $\meas(B^n) \ge \tfrac{\pi}{9} - \tfrac{\epsilon_n}{2} > 0$. Using the fact that $13$ is an escape radius for our sequence $\Pm$ gives that each $B^n$ is contained in the iterated filled Julia set $\K_{M_{n-1} + s_n}$. 

Now let $z_0 \in B^1$ and let $z^{n_k}$ be a sequence in our dense subset $(\Q \times \Q) \cup D^1$ which tends to $z_0$. By part 4. of the induction, this sequence lies in a union of the sets $I_i$. However, as the sets $B^{1,i}$ are open, it cannot lie in any finite union of these sets $I_i$ since such a sequence could not then converge to $z_0$. It then follows from part 5. of the induction that $z_0$ can be approximated by points in ${\mathcal A}_{\infty, s_1}$ and thus must lie in the Julia set $\J_{s_1}$ whence $\J_{s_1}$ has positive area as desired. By Theorem 1.1, all the iterated Julia sets $\Jm$, $m \ge 0$ must then also have positive area.  

To construct our invariant sequence of line fields, we first observe that we cannot pull back a non-zero Beltrami differential at a critical value of a mapping and so we need to remove the grand orbits of any critical points from the sets where we define our line fields. To this end, let $\{p_{m_k}\}$ be a listing of all those critical points (if any) which lie on the corresponding iterated filled Julia set, i.e. $p_{m_k}$ is a critical point of $P_{m_k + 1}$ with $p_{m_k} \in \K_{m_k}$. For each $n \ge 1$, set 

\[ \tilde B^n = B^n \setminus \bigcup_k \G_{M_{n-1} + s_n}(p_{m_k})\]

and for each $m \ge 0$ set 

\[ \Jmt = \Jm \setminus \bigcup_k \G_m(p_{m_k}), \qquad \Kmt = \Km \setminus \bigcup_k \G_m(p_{m_k}).\]

For each $n$, let $\tilde F^n$ be the set $Q_{M_{n-1} + s_n}^{\circ -1}(\tilde B^n) \subset E^n$. By complete invariance $\tilde F^n \subset \K$ for each $n$ and by the first, second  and third induction hypotheses, we have that $\meas(\K \setminus \tilde F^n) < 2^{-n}$. By part 7. of the induction, the sets $\tilde F^n$ are increasing in $n$ and if we then set $\tilde F = \cup_{n \ge 1}\tilde F^n$, then $\tilde F$ is a subset of $\K$ with $\meas(\tilde F) = \meas(\K)$.
 
Now define a line field $\tilde \mu_{s_1}$ on the set $\tilde B^1$ at time $s_1$, e.g. by setting $\tilde \mu_{s_1} = \dzbardz$ on $B^1$. By part 8. of our induction, the compositions $Q_{s_1, m}$, $m > s_1$ are all injective and so we can push $\tilde \mu_{s_1}$ forward to define a line field on all the sets $Q_{s_1, m}(\tilde B^1)$. On the other hand, for each $0 \le m < s_1$, we can use the composition $Q_{m,s_1}$ to define $\tilde \mu_m$ on the preimage of $\tilde B^1$ under this polynomial (note that from above we have avoided any potential problems with pulling back at critical values of $Q_{m, s_1}$). In particular, we can define $\tilde \mu_0$ on $\tilde F^1$. 

Now suppose that we have defined $\tilde \mu_{M_{n-1} + s_n}$ on the set $\tilde B^n$. Suppose also that, using part 8. of the induction again, by pushing forward, we have defined $\tilde \mu_m$  on the sets $Q_{M_{n-1} + s_n, m}(\tilde B^n)$ for $m > M_{n-1} + s_n$. Lastly, suppose that by pulling back, we have defined $\tilde \mu_m$ on the inverse images $Q_{m,M_{n-1} + m_m}^{\circ -1}(\tilde B^n)$ for all $0 \le m < M_{n-1} + s_n$, including the set $\tilde F^n$ at time $0$. 

At time $M_n + s_{n+1}$, $\tilde \mu_m$ is defined on $Q_{M_{n-1} + s_n, M_n + s_{n+1}}(\tilde B^n)$ which, by part 7. of the induction, is a subset of $\tilde B^{n+1}$. We then set $\tilde \mu_{M_n + s_{n+1}}$ to be simply $\dzbardz$ on the remainder of $\tilde B^{n+1}$, use part 8. of the induction to define $\tilde \mu_m$ on the forward images of $\tilde B^{n+1}$ at all times $m > M_n + s_{n+1}$ and then pull back as before to define $\tilde \mu_m$ on the preimages of this set at all previous times $0 \le m \le M_{n-1} + s_n$. Note that this will not change the existing definition of $\tilde \mu_m$ on the set $Q_{M_{n-1} + s_n, m}(\tilde B^n)$ for $m \ge M_{n-1} + s_n$ or on the set $Q_{m,M_{n-1}+s_n}^{\circ -1}(\tilde B^n)$ for $0 \le m < M_{n-1} + s_n$.

Proceeding in this way, we inductively define a sequence of line fields $\{\tilde \mu_m\}_{m=0}^\infty$. Our initial line field $\tilde \mu_0$ is then defined on $\tilde F$ and for each $m \ge 1$, $\tilde \mu_m$ is defined on $\tilde F_m := Q_m(\tilde F)$ which is a subset of the iterated filled Julia set $\Km$ of full measure by Theorem 1.1.

By construction, this sequence of line fields is invariant and it follows that on the set $\tilde G_m = \Kmt \setminus \tilde F_m$, if $z \in \tilde G_m$, then ${\mathcal G}_m(z) \subset \tilde G_m$. We may then pick a point $z$ of each grand orbit in $\tilde G_0$ on which we define $\tilde \mu_0$ to be $\dzbardz$ and extend it to the entire grand orbit of $\Gm(z)$ at every time $m \ge 0$ by pushing forwards and then pulling back. By the completeness of Lebesgue measure, the resulting extension of $\{\tilde \mu_m\}_{m=0}^\infty$ will still be measurable and each $\tilde \mu_m$ is now defined on the whole of $\Kmt$ which is the filled iterated Julia set $\K_m$ at most with countably many exceptions arising from the grand orbits of critical points. 

The very last step is to define an invariant sequence of line fields which is supported on subsets of the Julia rather than the filled Julia sets. We make the obvious restriction where for each $m \ge 0$, we define $\mu_m$ by 

\[  \mu_m (z) = \left \{ \begin{array}{r@{\: ,\quad}l} 
\tilde \mu_m(z) & z \in \Jmt \\
0 & z \notin \Jmt \\
\end{array} \right .\]

As all the iterated Julia sets have positive area and the grand orbits of any critical points on the Julia set give us at most a countably infinite set, we obtain an invariant sequence of measurable line fields whose supports are as desired. $\Box$

We conclude with a few remarks. Firstly, as our example relies so heavily on parabolic dynamics, the behaviour of the iterates on Julia sets is far from being hyperbolic and any suitable non-autonomous analogue of the radial Julia set is likely to be very small. This is of interest in view of Theorem 3.17 in \cite{McM} which implies that for a rational function with an invariant line field supported on its Julia set, almost every point in the Julia set approaches the postcritical set arbitrarily closely (Rempe and Van Strien have results in a similar direction in \cite{RVS}). 

It seems reasonable that we could perturb our sequence so as to ensure that there are no critical points on the iterated Julia sets, which would make for a slightly cleaner statement. Essentially what is required is estimates which state that the inverse images of sets of small measure still have small measure not just for one sequence of polynomials, but for all polynomial sequences with uniform bounds. Similar estimates are also needed for saying that suitable inverse images of points which are close remain close. This would then allow us to perturb our sequence slightly without destroying the required properties for the line field to exist. 

The polynomial $P^1 = \tfrac{1}{2}z(z + 1)$ was chosen simply because it had an attracting fixed point which was not superattracting. Any quadratic polynomial which possesses this property would also work and it is not too difficult to see from this that we can in fact choose our sequence $\Pm$ to be an arbitrary small perturbation of the constant sequence given by $P^2 = z(z + 1)$. Lastly, in view of Theorem 2.1 in \cite{Com3}, we may in fact assume the sequence in Theorem 1.2 is of the form $\{z^2 + c_m\}_{m=1}^\infty$.

Although we constructed an invariant sequence of line fields, we could just as easily have chosen Beltrami differentials $\mu_m(z) \dzbardz$ for which the $\linfty$ norms $|| \mu_m(z)||$ of the corresponding dilatations had absolute value uniformly less than $1$. Following the example of polynomial-like mappings, it seems probable that this would allow us to solve the Beltrami equation and conjugate $\Pm$ in the appropriate sense to another polynomial sequence $\Pmt$. Further, if the norms $|| \mu_m(z)||$ are uniformly small, then the sequence $\Pmt$ should be close to $\Pm$. The question is, then, how much of the nearby sequences to $\Pm$ can be obtained in this way? One has great freedom in the exact choice of invariant ellipse fields and if we could realise all sequences in a parameter neighbourhood of $\Pm$ in this way, we would then have a counterexample to the conjecture of density of hyperbolicity for non-autnonomous iteration, at least for sequences of quadratic polynomials. However, given the delicate nature of our construction and especially the heavy reliance on the parabolic dynamics associated with $P^2$, this seems rather unlikely.


\begin{thebibliography}{99}

\bibitem{BB} Rainer Br\"uck, Matthias B\"uger, \textsl{Generalized iteration}, 
Computational Methods and Function Theory Volume 3 (2003), {\bf 1}, 201-252.

\bibitem{BBR} Rainer Br\"uck, Matthias B\"uger, Stefan Reitz, \textsl{Random
iterations of\\ polynomials of the form $z^2 + c_n$: connectedness of Julia
sets}, Ergodic Theory Dynamical Systems (1999), {\bf 19}, 1221-1231.

\bibitem{CG} L. Carleson, T. W. Gamelin, ``Complex 
Dynamics'', Springer Verlag,\\ Universitext: Tracts in Mathematics (1993). 

\bibitem{Com1} M. Comerford, \textsl{Properties of Julia sets for the 
arbitrary composition of monic polynomials with uniformly bounded 
coefficients} Ph.D. Thesis, Yale University 2001.

\bibitem{Com2} M. Comerford, \textsl{A survey of results in random iteration} Proceedings Symposia in Pure Mathematics, American Mathematical Society, 2004. 

\bibitem{Com3} M. Comerford, \textsl{Conjugacy and counterexample in random iteration} Pac. J. of Math. {\bf 211} (2003), 69-80.

\bibitem{EL} Alexandre \`Eremenko and Mikhail Yu. Lyubich, \textsl{Examples of entire functions with pathological dynamics} J. London Math. Soc. (2) {\bf 36} (1987), no. 3, 458-468.

\bibitem{FS} J. E. Fornaess, Nessim Sibony, \textsl{Random iterations 
of rational functions} Ergodic Theory Dynamical Systems {\bf 11} (1991), 687-708.

\bibitem{Mil} J. Milnor, ``Dynamics in one Complex Variable'', 
Institute for\\ Mathematical Sciences, SUNY, Stony Brook, N.Y. 1991. 

\bibitem{McM} Curtis T. McMullen, ``Complex Dynamics and Renormalization'',\\
Annals of Mathematics Study 135, Princeton University Press, 1994.

\bibitem{MSS} R. Ma\~n\'e, P. Sad, D. Sullivan, \textsl{On the dynamics of rational maps}\\ Ann. Sc. de l'Ecole Normale Sup\'erieure {\bf 16}(1983), 193-217. 

\bibitem{RVS} L. Rempe, S. Van Strien, \textsl{Absence of line fields and Ma\~n\'e's theorem for nonrecurrent transcendental functions,} Transactions of the American Mathematical Society, {\bf 363} (2011), 203-228.

\bibitem{Wang} Xiaoguang Wang, \textsl{Rational Maps Admitting Meromorphic Invariant Line Fields},  Bull. Aust. Math. Soc. {\bf 80} (2009), 454Ð461.


\end{thebibliography}
\end{document}